\documentclass[graybox,natbib]{svmult}

\usepackage{mathptmx}       
\usepackage{helvet}         
\usepackage{courier}        
\usepackage{makeidx}         
\usepackage{graphicx}        
\usepackage{multicol}        
\usepackage[bottom]{footmisc}
\usepackage{amstext}
\usepackage{amsmath}
\usepackage[utf8]{inputenc}     
\usepackage{booktabs}
\usepackage{url}                
\setcitestyle{round,longnamesfirst}

\newcommand{\discuss}[1]{\begin{svgraybox}\textbf{Discussion Point:} #1\end{svgraybox}}

\definecolor{helena}{rgb}{.2,.8,.4}

\makeindex             

\begin{document}

\title*{Nonnegative Factorization of a Data Matrix as a Motivational Example for Basic Linear Algebra}
\titlerunning{Nonnegative Factorization of a Data Matrix as a Motivational Example \ldots}
\author{Barak A. Pearlmutter \and Helena Šmigoc}
\authorrunning{B. A. Pearlmutter and H. Šmigoc}
\institute{Prof.\ Barak A. Pearlmutter \at Department of Computer Science, Maynooth University, Ireland \newline \email{barak@pearlmutter.net}
\and Dr.\ Helena Šmigoc \at School of Mathematics and Statistics, UCD Dublin, Ireland \newline \email{helena.smigoc@ucd.ie}}
%
%
\maketitle

\abstract{ We present a motivating example for matrix multiplication based on factoring a data matrix. Traditionally, matrix multiplication is motivated by applications in physics: composing rigid transformations, scaling, sheering, etc. We present an engaging modern example which naturally motivates a variety of matrix manipulations, and a variety of different ways of viewing matrix multiplication. We exhibit a low-rank non-negative decomposition (NMF) of a ``data matrix'' whose entries are word frequencies across a corpus of documents. We then explore the meaning of the entries in the decomposition, find natural interpretations of intermediate quantities that arise in several different ways of writing the matrix product, and show the utility of various matrix operations. This example gives the students a glimpse of the power of an advanced linear algebraic technique used in modern data science.}

\keywords{Nonnegative Matrix Factorization (NMF);
Topic Modeling;
Data Mining;
Matrix Multiplication}

\section{Introduction}
\label{sec:intro}

Examples are an essential part of teaching any mathematical subject. They serve a range of purposes, from checking understanding and deepening knowledge to giving a broader view of the subject and its applications. There are an abundance of examples available in the literature, covering every topic of any basic linear algebra course. However, it is not so easy to find examples that give an insight into the current development of the subject and are at the same time accessible to students. As the applications of linear algebra are rapidly expanding, and several new developments in the subject are motivated by applications, examples showcasing current applications of the subject are of particular interest.

Because of its utility in other domains, linear algebra is a classical subject routinely taught to students not majoring in mathematics. It is a prerequisite not just for advanced mathematics but also for undergraduate degrees in Engineering, Physics, Computer Science, Biology, Chemistry,  Business, Statistics, and the like. Those students in particular benefit from learning from examples, and appreciate seeing interesting applications of the material they are learning. The benefits of using models to introduce mathematical concepts has been studied \citep{Lesh2005}, and models focusing on different concepts from linear algebra are available \citep{Possani20102125, Salgado2015100, Trigueros20131779}.

While very simple examples are essential when introducing a topic, examples of applications presented in classrooms often seem contrived. For example, students' knowledge of economics and agriculture is sufficiently sophisticated that simple linear examples of acreage under cultivation invite criticism.
On the other hand, it is impossible to bring to the classroom, for instance, deep applications of linear algebra in genetics \citep{Ponnapalli-etal-2011a}, since most likely neither the instructor nor the students have the necessary background to really understand how they work. To quote \citet{Stewart03},
\begin{quotation}
  While it is true that linear algebra can simplify the solution to many problems, this is only true for those who are very familiar with the subject area. In contrast, the first year university student has a long way to go before being able to see the whole picture.
\end{quotation}

The press is full of stories about data science: analysis of large corpora of data. Some of these lend themselves naturally to use as motivating examples for various concepts in linear algebra. For example, the Netflix challenge can be viewed as a problem in matrix completion, where a company was highly motivated to recover a low rank decomposition of an almost entirely incomplete matrix of movie ratings.

We present less abstract example, in which matrix multiplication is explicated by examination of a nonnegative decomposition of a term-by-document matrix.  This particular example vividly illustrates  various views of matrix multiplication (as composition of linear functions; as a sum of outer products of columns with rows; and as a table of inner products of rows with columns), while using only primitive concepts.  It also previews and motivates a variety of more advanced concepts (the general algebraic concept of factoring, the notion of rank, approximation and norms, iterative numeric algorithms, constraints like element-wise non-negativity, and column-stochastic matrices), helping sketch the outlines of richer material covered in more advanced courses.

Although intuitive and implemented by a very short algorithm, the technique discussed (NMF) is far from a toy: it has enjoyed a myriad of accessible and engaging applications \citep{ASARI-ETAL-2006A, Cancer-NMF, ECG-EMG-NMF, OGRADY-PEARLMUTTER-2008a, HIV-NMF, Smaragdis-Brown-2003a, Wilson-etal-2008a}.  For this reason, the example we present serves to give a taste of an interesting and accessible application of linear algebra.  Although briefly presented in this document for the sake of completeness, we do not suggest attempting to derive the method in the classroom, leaving that too as motivation for the pursuit of more advanced study.


\smallskip

The mathematical notation used below is standard. For example, $e_i$ denotes the vector of appropriate size with $i$-th entry equal to one and other entries equal to zero.
The ``discussion point'' boxes are intended to be illustrative, and can be used for classroom discussion, project-based learning, or as the basis for assignments.


\section{Term-by-Document Matrix: A Small Example}
\label{sec:smallexample}

Numeric data organised in a tabular format is something we are all familiar with in our daily lives. Everyone can understand a spreadsheet whose rows are indexed by products, columns by month, and whose entries contain sales. These are the matrices that students entering a linear algebra course have already seen. In data science, tabular data of this sort is known as a ``data matrix''.

A data matrix of interest in library science is a tabulation of word frequencies by documents. Rows are indexed by words, columns by documents, and the entries of a matrix are the number of times a given word appears in a given document. This particular kind of data matrix is sometimes called a term-by-document matrix. Although this matrix completely ignores the actual arrangement of words within each document (i.e., it is a bag-of-words model), it still contains sufficient information to allow interesting structure to be discovered.

There are several ways in which matrices and matrix multiplication can be introduced in the classroom. Term-by-document matrices can be one of the examples given to the class, starting with a small example that can be given on a board. 
In the classroom we can show a pre-prepared example, which can be built on by an assignment in which students have the freedom to chose the documents they want to consider. Since the search function in browsers automatically counts the number of times a word appears on a page, such an assignment is not necessarily time demanding.

Here we present an example where the documents are the Wikipedia entries for the four most venomous animals in the world
(\emph{Box Jellyfish}\footnote{\url{https://en.wikipedia.org/wiki/Box_jellyfish}},
\emph{King Cobra}\footnote{\url{https://en.wikipedia.org/wiki/King_cobra}},
\emph{Marbled Cone Snail}\footnote{\url{https://en.wikipedia.org/wiki/Conus_marmoreus}},
\emph{Blue-Ringed Octopus}\footnote{\url{https://en.wikipedia.org/wiki/Blue-ringed_octopus}})
and we consider only five terms
(\emph{venom},
\emph{death},
\emph{danger},
\emph{survive},
\emph{Madagascar}).
This gives us Table~\ref{tab:Toy Example}.
\begin{table}
\caption{Term-by-Document Matrix of the Four Most Venomous Animals}
\label{tab:Toy Example}  
\renewcommand{\arraystretch}{1.2}
\renewcommand{\tabcolsep}{1em}
\begin{tabular}{|l r | r@{\hspace{2em}} r@{\hspace{1.4em}} r@{\hspace{1.5em}} r@{\hspace{2.0em}} |}
  \cline{3-6}
\multicolumn{2}{c|}{} &  \multicolumn{4}{c|}{\textbf{Documents}} \\
\multicolumn{2}{c|}{}
& \multicolumn{1}{c}{Jellyfish}
& \multicolumn{1}{c}{Cobra}
& \multicolumn{1}{c}{Snail}
& \multicolumn{1}{c|}{Octopus} \\
\hline
\textbf{Terms}
& venom & 32 & 44 & 1 & 18 \\
& death & 9 & 3 & 0 & 2 \\
& danger & 6 & 4 & 0 & 4\\
& survive & 2 & 0 & 0 & 1\\
& Madagascar & 0 & 0 & 2 & 0\\
\hline
\end{tabular}
\end{table}

Going from the table to the matrix
\[
A=\left(
\begin{array}{cccc}
 32 & 44 & 1 & 18 \\
 9 & 3 & 0 & 2 \\
 6 & 4 & 0 & 4 \\
 2 & 0 & 0 & 1 \\
 0 & 0 & 2 & 0 \\
\end{array}
\right)
\]
we can lead the discussion in several directions. A representative set of questions is given below. The questions are of course trivial to answer without referring to matrices.  The simplicity of the questions makes it easy for students to understand the corresponding matrix operations and motivates them to think about extensions to more involved tasks.  

\discuss{Determine the frequency of terms appearing in the first document, in the third document, in the first or third document. What is the frequency of terms in all the documents together?}

The above questions can all be answered using multiplication of a matrix by a column vector.
\begin{align*}
A\left(\begin{array}{c}1 \\ 0 \\ 0 \\ 0 \end{array}\right)
&=\left(\begin{array}{c}32 \\ 9 \\ 6 \\ 2 \\ 0 \end{array}\right)
&
A\left(\begin{array}{c}0 \\ 0 \\ 1 \\ 0 \end{array}\right)
&=\left(\begin{array}{c}1 \\ 0 \\ 0 \\ 0 \\ 2\end{array}\right)
\\[2ex]
A\left(\begin{array}{c}1 \\ 0 \\ 1 \\ 0 \end{array}\right)
&=\left(\begin{array}{c}33 \\ 9 \\ 6 \\ 2 \\ 2 \end{array}\right)
&
A\left(\begin{array}{c}1 \\ 1 \\ 1 \\ 1 \end{array}\right)
&=\left(\begin{array}{c}95 \\ 14 \\ 14 \\ 3 \\ 2 \end{array}\right)
\end{align*}
On this small example matrix multiplication, while illustrative, does not help with efficiency of obtaining an answer. However, we can lead the students to think further.   

\discuss{Can you think about other questions about the set of documents that can be answered using matrix multiplication? (E.g., differences in word frequencies.)
How would one extract information from very large datasets? (This is for the computer science students in the class: strategies for assembling, representing, storing, and operating upon a very large data matrix.)}

At this point the students can appreciate that in order to get the information\footnote{The ``information'' here can be viewed as histograms over either documents or terms, which is something the students should be comfortable with.} about the terms in the $i$-th document we need to multiply $A$ by $e_i$, to find the information about the terms in documents $i$, $j$ and $k$ we need to multiply $A$ by $e_{i}+e_{j}+e_{k}$, or equivalently, add $Ae_i$, $Ae_j$, and $Ae_k$. We can view the matrix $A$ as a transformation that takes information about documents ($i$, $j$, $k$) to information about terms ($Ae_i$, $Ae_j$, $Ae_k$).
\[
 \text{words} \stackrel{~A}{\longleftarrow} \text{documents}
\]
Furthermore, we can remark that this transformation obeys certain rules
\[
A(e_i+e_j+e_k)=Ae_i+Ae_j+Ae_k
\]
which can be developed into the definition of linearity.
We continue the discussion by presenting the transpose matrix.
\[
A^T=\left(
\begin{array}{ccccc}
 32 & 9 & 6 & 2 & 0 \\
 44 & 3 & 4 & 0 & 0 \\
 1 & 0 & 0 & 0 & 2 \\
 18 & 2 & 4 & 1 & 0 \\
\end{array}
\right)
\]

\discuss{Which documents contain the third term, the fifth term, the third or the fifth term?} 
\begin{align*}
A^T\left(\begin{array}{c}0 \\ 0 \\ 1 \\ 0 \\ 0 \end{array}\right)
&=\left(
\begin{array}{c}
6 \\ 4 \\ 0 \\ 4 \\
\end{array}
\right)
&
A^T\left(\begin{array}{c}0 \\ 0 \\ 0 \\ 0 \\ 1 \end{array}\right)
&=\left(
\begin{array}{c}
 0 \\ 0 \\ 2 \\ 0 \\
\end{array}
\right)
&
A^T\left(\begin{array}{c}0 \\ 0 \\ 1 \\ 0 \\ 1 \end{array}\right)
&=\left(
\begin{array}{c}
 6 \\ 4 \\ 2 \\ 4  \\
\end{array}
\right)
\end{align*}
Students can see the similarity with the discussion above. To find out how the $i$-th term is featured in documents we need to multiply $e_i$ by $A^T.$
\[
 \text{documents} \stackrel{\;\;\;\;A^T}{\longleftarrow} \text{words}
\]
Further questions  can be discussed in this framework, touching upon elementary ideas not routinely discussed in the first course on linear algebra, such as non-negativity and sparsity.\footnote{Sparsity is of particular importance in computer science, where it impacts the representation and manipulation of both matrices and graphs.}

\discuss{Note that the term Madagascar only appears in the third document. Can we draw any conclusions from this?}

\discuss{If we were to make a table that includes all the terms that appear in at least one of the four documents, would we expect most of the entries in the matrix to be equal to zero? Why?}

\discuss{Note that all the elements in the matrix are non-negative integers. Can you think of any other tables with only non-negative integers? How about tables containing only non-negative real elements? Can you think about any other conditions on the entries that are imposed naturally in a particular setting?}


\discuss{In class, we usually label the rows and columns of a matrix with successive integers: $1, \ldots, n$.
These are generally used as ``nominal numbers'', meaning only their identities are important---like building numbers, course numbers, or social security numbers.
And when we write $\sum_{i=1}^n$, what we usually mean is really $\sum_{i \in \text{rows}}$.
We can change most of our formulas to use this convention.
But in actual applications, as in the example here, often the rows and columns have natural labels: names of chemicals, words, documents, people, months, cities, \emph{etc.}
When this holds, we can use these labels instead of numbers as indices.
And we can freely rearrange the rows and columns, keeping their labels, while still representing the same underlying mathematical object: the same matrix.}
\noindent
This point is illustrated by a term-by-document matrix, which has rows labeled by terms and columns labeled by documents. Let us look at another example. The two tables below contain movie ratings given by four users to five movies: 
\newcommand{\personA}[0]{Cindy}
\newcommand{\personB}[0]{Dora}
\newcommand{\personC}[0]{Alice}
\newcommand{\personD}[0]{Becky}
\newcommand{\movieA}[0]{\emph{Alien}}
\newcommand{\movieB}[0]{\emph{Jaws}}
\newcommand{\movieC}[0]{\emph{Beetlejuice}}
\newcommand{\movieD}[0]{\emph{Animal House}}
\newcommand{\movieE}[0]{\emph{Life of Brian}}
\begin{table}
\caption{Labeling Rows and Columns}
\label{tab:Movie Ratings}
\begin{tabular}{r|cccc|}
\multicolumn{1}{r}{}
        & \personC      & \personD      & \personA      &
\multicolumn{1}{r}{
                                                          \personB}\\\cline{2-5}
\movieA &       4       &        1      &       4       &       5  \\
\movieD &       1       &        5      &       4       &       2  \\
\movieC &       2       &        2      &       5       &       3  \\
\movieB &       5       &        1      &       5       &       5  \\
\movieE &       1       &        5      &       5       &       1  \\
\cline{2-5}
\end{tabular}
\hfill
\begin{tabular}{r|cccc|}
\multicolumn{1}{r}{}
        & \personD      & \personA      & \personB      &
\multicolumn{1}{r}{
                                                          \personC}\\\cline{2-5}
\movieD &          5      &    4       &        2      &       1   \\
\movieE &          5      &    5       &        1      &       1   \\
\movieC &          2      &    5       &        3      &       2   \\
\movieB &          1      &    5       &        5      &       5   \\
\movieA &          1      &    4       &        5      &       4   \\
\cline{2-5}
\end{tabular}
\end{table}

\discuss{Compare the two tables. Do they contain the same information? Can you figure out the principle behind the ordering of rows and columns on the left and on the right?}

\section{Matrix Factorization}

Students are familiar with the idea of factoring an integer as a product of prime numbers. Writing $6=3\times 2$ gives us some information about the number $6$. Another example is factoring a polynomial. Writing
\(
x^4-10 x^3+35 x^2-50 x+24
\)
as 
\(
(x-1)(x-2)(x-3)(x-4)
\)
uncovers useful information. Both prime factor decomposition and factoring a polynomial are in general hard to do. Given two integers it is straightforward to find their product, but there is no known efficient algorithm for integer factorization.

This concept can, in some sense, be extended to matrices. Given a matrix, we want to write it as a product of two (or more) ``simpler'' matrices. There are several ways this can be done. A wide range of factorizations of matrices are used in applications, where---depending on the application---different properties of the factors are desired. An example that can be presented in the classroom is given below. The matrix
$$
A=\left(
\begin{array}{cccc}
 2 & -1 & 1 & 2 \\
 -1 & 1 & -2 & -1 \\
 1 & -2 & 5 & 1 \\
 2 & -1 & 1 & 2 \\
\end{array}
\right)$$
can be factored in several ways: 
\begin{align*}
A&=\frac{1}{21} \left(
\begin{array}{cccc}
 1 & 1 & -1 & 1 \\
 -1 & 0 & 0 & 3 \\
 2 & -1 & 0 & 1 \\
 1 & 1 & 1 & 0 \\
\end{array}
\right)
\left(
\begin{array}{cccc}
 7 & 0 & 0 & 0 \\
 0 & 3 & 0 & 0 \\
 0 & 0 & 0 & 0 \\
 0 & 0 & 0 & 0 \\
\end{array}
\right)
\left(
\begin{array}{cccc}
 3 & -3 & 6 & 3 \\
 7 & 0 & -7 & 7 \\
 -10 & 3 & 1 & 11 \\
 1 & 6 & 2 & 1 \\
\end{array}
\right) \\
&=\left(
\begin{array}{cc}
 1 & 1 \\
 -1 & 0 \\
 2 & -1 \\
 1 & 1 \\
\end{array}
\right)
\left(
\begin{array}{cccc}
 1 & -1 & 2 & 1 \\
 1 & 0 & -1 & 1 \\
\end{array}
\right) \\
&=\left(
\begin{array}{cccc}
 2 & -1 & 0 & 0 \\
 -1 & 1 & 0 & 0 \\
 1 & -2 & 1 & 0 \\
 2 & -1 & 0 & 1 \\
\end{array}
\right)
\left(
\begin{array}{cccc}
 1 & 0 & -1 & 1 \\
 0 & 1 & -3 & 0 \\
 0 & 0 & 0 & 0 \\
 0 & 0 & 0 & 0 \\
\end{array}
\right)
\end{align*}
The students may not have the mathematical tools to develop the factorizations above, but we can ask them to check their correctness, and to explore properties of the factors. 

Demands from applications frequently impose conditions on the factors that are too strong to be satisfied exactly. For example, not every matrix can be written as a product of a column by a row. Or more generally, not every matrix can be written as a product of two low rank matrices. If we are unwilling to relax the conditions, we need to resort to approximate factorizations. Let us consider the matrix 
$$A=\left(
\begin{array}{cccc}
 1 & 1 & 1 & 1.01 \\
 1 & 1 & 1.01 & 1 \\
 1 & 1.01 & 1 & 1 \\
 1.01 & 1 & 1 & 1 \\
\end{array}
\right).$$
Using elementary tools one can check that $A$ cannot be written as a product of a column and a row. 

\discuss{Can we find a matrix that is close to the matrix $A$ that can be written as a product of a column by a row?}

\noindent Students are likely to come up with the following solution: 
$$A_1=\left(
\begin{array}{cccc}
 1 & 1 & 1 & 1 \\
 1 & 1 & 1 & 1 \\
 1 & 1 & 1 & 1 \\
 1 & 1 & 1 & 1 \\
\end{array}
\right)=\left(
\begin{array}{c}
 1 \\
 1 \\
 1 \\
 1 \\
\end{array}
\right)\left(
\begin{array}{cccc}
 1 & 1 & 1 & 1 \\
\end{array}
\right),$$ 
and we may present another one: 
$$A_2=\left(
\begin{array}{cccc}
 1.0025 & 1.0025 & 1.0025 & 1.0025 \\
 1.0025 & 1.0025 & 1.0025 & 1.0025 \\
 1.0025 & 1.0025 & 1.0025 & 1.0025 \\
 1.0025 & 1.0025 & 1.0025 & 1.0025 \\
\end{array}
\right)=\left(
\begin{array}{c}
 1 \\
 1 \\
 1 \\
 1 \\
\end{array}
\right)\left(
\begin{array}{cccc}
 1.0025 & 1.0025 & 1.0025 & 1.0025 \\
\end{array}
\right).$$

\discuss{Which solutions is better? What does it mean for a matrix $B$ to be close to $A$?}

\noindent This question could be an introduction to the concept of matrix norms.

Factorizations of matrices have been developed that are used in applications, with the aim of uncovering hidden structure. We present a factorization that requires the factors to be non-negative and of a given low rank $r$. Those conditions are too strong, so the factorization will not be exact. This means that given a matrix $V$, we obtain matrices $W$ and $H$ so that the matrix $\hat{V}=WH$ is in some sense close to $V$. 
\[
  \fbox{\rule{3em}{0ex}\rule{0em}{5ex}\raisebox{1.7ex}{$V$}\rule{3em}{0ex}} \; \raisebox{1.7ex}{$\approx$} \; \fbox{\rule{0em}{5ex}\raisebox{1.7ex}{$W$}} \;\fbox{\rule{0em}{1.8ex}\rule{3em}{0ex}\raisebox{0.15ex}{$H$}\rule{3em}{0ex}}
\]

Non-negative matrix factorization, or NMF \citep{PAATERO-TAPPER-1994, LEE-SEUNG-99A, Wang-Zhang-2013a}, is a class of techniques for \emph{approximately} factoring a matrix of non-negative numbers into the product of two such matrices: given an entry-wise non-negative $n \times m$ matrix $V$, find two entry-wise non-negative matrices $W$ and $H$, of sizes $n \times r$ and $r \times m$, such that $V \approx W H$.
(Even after a value for $r$ has been chosen, and an appropriate measure of similarity of two matrices has been chosen, there can be many possible solutions. However, popular NMF algorithms empirically usually find good solutions, a phenomenon which has been the subject of considerable analysis \citep{DONOHO2003}.)

Let us look at the nonnegative matrix factorization of the matrix that corresponds to the left side of Table~\ref{tab:Movie Ratings}:
$$A=\left(\begin{array}{cccc} 4  & 5  & 4  & 1 \\ 5  & 5  & 5  & 1 \\ 5  & 3  & 2  & 2 \\ 4  & 2  & 1  & 5 \\ 5  & 1  & 1  & 5  \end{array}\right)$$ 
First we take $r$ to be equal to one. That means that we want to approximate $A$ by a product of a nonnegative column and a nonnegative row. The algorithm returns the following result: 
\[
W_1=\left(\begin{array}{c} 7.137\\ 8.214\\ 6.398\\ 5.974\\ 6.155 \end{array}\right)
\hspace{4em}
H_1=\left(\begin{array}{cccc} 0.6709 & 0.4898 & 0.406 & 0.381 \end{array}\right)
\]
\[
A-W_1H_1=\left(\begin{array}{cccc} -0.7885 & 1.504 & 1.102 & -1.72\\ -0.511 & 0.977 & 1.665 & -2.13\\ 0.7077 & -0.1334 & -0.5976 & -0.4378\\ -0.008175 & -0.926 & -1.426 & 2.724\\ 0.8703 & -2.015 & -1.499 & 2.655 \end{array}\right)
\]
Taking $r=2$ we get:
\[
W_2=\left(\begin{array}{cc} 6.968 & 1.086\\ 7.908 & 1.364\\ 3.763 & 3.558\\ 0.5117 & 6.448\\ 0 & 7.197 \end{array}\right)
\hspace{4em}
H_2=\left(\begin{array}{cccc} 0.5171 & 0.6379 & 0.5707 & 0\\ 0.6658 & 0.1897 & 0.1095 & 0.7133 \end{array}\right)
\]
\[
A-W_2H_2=\left(\begin{array}{cccc} -0.3256 & 0.3496 & -0.09564 & 0.2257\\ 0.002296 & -0.3033 & 0.337 & 0.02677\\ 0.6857 & -0.07508 & -0.5374 & -0.5376\\ -0.5576 & 0.4506 & 0.001623 & 0.4004\\ 0.2088 & -0.365 & 0.2117 & -0.1333 \end{array}\right)
\]

\discuss{Compare $A-W_1H_1$ and $A-W_2H_2$. Can you find a nonnegative factorization of $A$ for $r=4$?}

Let us have a closer look at $W_2$ and $H_2$. Recall that the rows of $W_2$ correspond to movies, and the columns of $H_2$ correspond to users.
\discuss{Can we give sensible labels to the columns of $W_2$, or equivalently, the rows of $H_2$.}
\noindent
Note that the highest values in the first column of $W_2$ correspond to movies \movieA\ and \movieB, while the highest values in the second column correspond to movies \movieD\ and \movieE.  Based on this, we may agree to label the first column ``Horror'', and the second column ``Comedy''.  The rows of $H_2$ are labeled correspondingly. Values in $H_2$ can now be interpreted in the following way. \personA\ likes both horror and comedy movies, \personB\ and \personC\ prefer horror movies, and \personD\ likes comedy, but not horror movies.  Matrix factorization uncovered genre for our movies. 

In the context of the small example, we can look at $V$ as a transformation from ``movies'' to ``people''. Now we have the third notion appearing: ``genres''. The matrix $H$ can be seen as a transformation that takes movies to genres, and $W$ takes genres to people.
\[
 \text{people}
 \overbrace{\stackrel{~W}{\longleftarrow} \text{genres}\stackrel{~H}{\longleftarrow}}^{V}
 \text{movies}
\]
We may remark to the students that they are justified in finding this example a bit contrived. The is example is too small (and also made up) to be very convincing. We give a larger example based on term-by-document matrix later in the chapter.  


\bigskip

The formula for matrix multiplication, $A=BC$,
\[
a_{ij}=\sum_{k=1}^m b_{ik}c_{kj}
\]
can be intimidating to students.

From the point of view of our example, where the rows of $V$ are indexed by ``movies'' and the columns by ``people'', we can write down the same formula in the following way: 
\[
 \hat{v}_{\text{movie},\text{person}}=\sum_{g\in\text{genre}}w_{\text{movie},g} \; h_{g, \text{person}}
\]
The entry in the matrix $\hat{V}$ that represents a rating of a chosen ``movie'' to a chosen ``person'' is computed by summing up the product of how much the person likes genre $g$ and how much the movie is in genre $g$, over all genres.
This process can be depicted graphically:
\[
\left[\begin{array}{ccc}
    \vspace{3ex}&\hspace*{7em}&\\
    &\hat{V}&\\
    \vspace{3ex}&\hspace{4em}{\bullet}&
  \end{array}\right]
=
\left[\begin{array}{ccc}
    \vspace{3ex}&\hspace*{1em}&\\
    &\hspace{1.5ex}W\hspace{1ex}&\\
{\makebox[0ex][l]{\rule{6ex}{0.5ex}}}\vspace{3ex}&&
  \end{array}\right]
\left[\begin{array}{ccc}
    \vspace{0ex}&\hspace*{7em}&\\
    &H&\\
    \vspace{0ex}&\hspace{4em}\parbox[b][0ex]{0em}{{\rule{0.5ex}{6ex}}}&
  \end{array}\right]
\]

%

On the other hand we may notice that $\hat V$ is the sum of rank one matrices, each of them giving the contribution of a particular genre. 
\[
\left[\begin{array}{ccc}
    \vspace{3ex}&\hspace*{7em}&\\
    &\hat{V}&\\
    \vspace{3ex}&\hspace{4em}&
  \end{array}\right]
{=}
\sum_g
\left[\begin{array}{@{}c@{}}
    \vspace{3ex}\\
    W_{\bullet,g}\\
    \vspace{3ex}
  \end{array}\right]
\left[\begin{array}{ccc}
    \hspace{4ex}&H_{g,\bullet}&\hspace{4ex}\\
  \end{array}\right].
\]
This is written as: 
\begin{eqnarray*}
V&=&\sum_{g=1}^{r} (\text{$g$-th column of $W$}) \cdot (\text{$g$-th row of $H$})\\
&=&\sum_{g \in \text{genres}} \text{Ratings Matrix for genre $g$}\\
&=&\text{sum of rank-one per-genre matrices}
\end{eqnarray*}
This interpretation reinforces the power of the nonnegative matrix factorization. From a bundle of documents, it singles out a particular genre in way that agrees with our intuition in a surprisingly strong way.

Let us go back to the example of movie ratings discussed earlier. We have approximated our matrix $A$ by $W_2H_2$. Below this product is written as a sum of two rank one matrices: 
\begin{align*}
W_2H_2&=\left(\begin{array}{cc} 6.968 & 1.086\\ 7.908 & 1.364\\ 3.763 & 3.558\\ 0.5117 & 6.448\\ 0 & 7.197 \end{array}\right)\left(\begin{array}{cccc} 0.5171 & 0.6379 & 0.5707 & 0\\ 0.6658 & 0.1897 & 0.1095 & 0.7133 \end{array}\right)\\
&=\left(\begin{array}{c} 6.968\\ 7.908 \\ 3.763 \\ 0.5117 \\ 0  \end{array}\right)\left(\begin{array}{cccc} 0.5171 & 0.6379 & 0.5707 & 0 \end{array}\right) \\&\qquad +
\left(\begin{array}{cc}  1.086\\ 1.364\\  3.558\\  6.448\\  7.197 \end{array}\right)\left(\begin{array}{cccc} 0.6658 & 0.1897 & 0.1095 & 0.7133 \end{array}\right)\\
&=\left(\begin{array}{cccc} 3.603 & 4.445 & 3.977 & 0\\ 4.089 & 5.045 & 4.514 & 0\\ 1.946 & 2.4 & 2.148 & 0\\ 0.2646 & 0.3264 & 0.292 & 0\\ 0 & 0 & 0 & 0 \end{array}\right)
+\left(\begin{array}{cccc} 0.7227 & 0.2059 & 0.1189 & 0.7743\\ 0.9084 & 0.2588 & 0.1495 & 0.9732\\ 2.368 & 0.6747 & 0.3897 & 2.538\\ 4.293 & 1.223 & 0.7064 & 4.6\\ 4.791 & 1.365 & 0.7883 & 5.133 \end{array}\right)
\end{align*}

\begin{svgraybox}
  \textbf{Supplementary Advanced Material}

  \medskip

  \noindent To minimize the Frobenius norm, $\lVert V - W H \rVert_F$, a surprisingly simple method is available. Two updates can be iterated
  \begin{eqnarray*}
    H &:=& H \odot (W^T V \div W^T W H)
    \\
    W &:=& W \odot (V H^T \div W H H^T)
  \end{eqnarray*}
  where $\odot$ denotes elementwise multiplication and $\div$ denotes elementwise division.
  Typically after each update of $W$ its columns are normalized to unit sum.
\end{svgraybox}



\noindent
An enormous number of variations and embellishments of the basic NMF algorithm have been developed, with applications ranging from astronomy to zoology.

\section{Example using Module Descriptors}

We give an example of factoring a data matrix involving a corpus of
documents: module descriptors for $62$ mathematics modules that were
taught in the School of Mathematics and Statistics, University College
Dublin (UCD) in 2015.  Module descriptors are relatively short
documents that give overviews of the courses. Here are two
representative examples of module descriptors.

\begin{quotation}
\fbox{\parbox{0.9\columnwidth}{
\begin{center}
\textbf{Numbers and Functions}
\end{center}
\raggedright
\noindent This module is an introduction to the joys and challenges of mathematical reasoning and mathematical problem-solving, organised primarily around the theme of properties of the whole numbers.
It begins with an introduction to some basic notions of mathematics and logic such as proof by contradiction and mathematical induction.
It introduces the language of sets and functions, including injective surjective and bijective maps and the related notions of left-, right- and 2-sided inverses. Equivalence relations, equivalence classes.
It covers basic important principles in combinatorics such as the Principle of Inclusion-Exclusion and the Pigeonhole Principle.
The greater part of the module is devoted to number theory: integers, greatest common divisors, prime numbers, Euclid's algorithm, the Fundamental Theorem of Arithmetic, congruences, Fermat's theorem, Euler's theorem, and arithmetic modulo a prime and applications. The module concludes with some topics from elementary coding theory / cryptography such as the RSA encryption system. 
}}
\end{quotation}

\begin{quotation}
\fbox{\parbox{0.9\columnwidth}{
\begin{center}
\textbf{Groups, Rings and Fields}
\end{center}
\raggedright
\noindent This course will be an introduction to group theory, ring theory and field theory. We will cover the following topics: definition and examples of groups, subgroups, cosets and Lagrange's Theorem, the order of an element of a group, normal subgroups and quotient groups, group homomorphisms and the homomorphism theorem, more isomorphism theorems, definitions of a commutative ring with unity, integral domains and fields, units, irreducibles and primes in a ring, ideals and quotient rings, prime and maximal ideals, ring homomorphisms and the homomorphism theorem, polynomial rings, the division algorithm, gcd for polynomials, irreducible polynomials and field extensions. Time permitting, we may cover the Sylow theorems, solvable groups and further examples of groups.
}}
\end{quotation}

This set was chosen because the example was developed for the first linear algebra classes at UCD. All the students in those classes were quite familiar with the chosen set of documents, which they need to navigate each semester when choosing and registering for their modules. The number of documents is large enough that we can make a case for needing a computer to help navigate them, but small enough that we can exert some manual control and can be familiar with the entire corpus.

All the words that appear in any of the $62$ documents were collected. So-called stop words (common words like ``and'', ``the'', and the like), and all the words that appeared fewer than four times, were removed. Words were also down-cased and stemmed, so for example the terms \emph{Eigenvalue} and \emph{eigenvalues} are deemed equivalent. This resulted in a $290 \times 62$ matrix $V$ of word counts. Below we show some of the data, starting with lists of the most and the least frequent words. 

\begin{table}
\caption{Most and Least Frequent Words}
\label{tab:Most Frequent Words}       
\begin{tabular}{lr}\toprule
  \multicolumn{2}{c}{\textbf{Most Frequent Words}}\\
  \midrule
\textbf{word} & \textbf{count} \\\midrule
function &	117 \\
theorem &	75 \\
linear &	72 \\
matrix &	66 \\
theory &	58\\
equation & 53 \\
mathematical & 	52 \\
mathematics &	52 \\
understand &	49 \\
science &	48 \\
problem &	44 \\
\bottomrule
\end{tabular}
\hfill
\begin{tabular}{lllll}\toprule
  \multicolumn{5}{c}{\textbf{Least Frequent Words}}\\
  \midrule
addition &	
advanced &	 
arguments &	 
background &	 
behaviour \\
classify &	 
column &	 
computation &	 
constrained &	 
construct \\	 
continuous &	 
definite &	 
depth &	 
described &	 
directional \\	 
double &	 
elimination &	 
engineering &	 
evaluate &	 
expressions \\	 	 
flow &	 
foundations &	 
general &	 
importance &	 
independence \\	 
induction &	 
integrate &	 
interpret &	 
introduces &	 
known \\	 
manipulate &	 
max &	 
maxima &	 
min &	 
minima \\	 
nash &	 
nullity &	 
numerical &	 
original &	 
possible \\	 
prime &	 
quadratic &	 
range &	 
related &	 
riemann \\	 
row &	 
sample &	 
search &	 
significant &	 
solid \\ 
special &	 
stock & 
sum &	 
sylow &	 
together \\	 
uncountable &	 
variety\\
\bottomrule
\end{tabular}
\end{table}

%

Using the standard Octave \citep{Octave-4} package for NMF, the entry-wise nonnegative matrix $V$ is factored as
\[
V\approx WH,
\]
where $W$ is an entry-wise non-negative $290 \times r$ matrix and $H$ is an entry-wise nonnegative $r \times 62$ matrix. We will see that factoring a matrix in this way reveals a particular structure of the matrix which reveals something about the content of the original documents. Different values for $r$ are used below, and we will see how the information that we obtain changes as we increase $r$.

For easier interpretation, the entries in each column of $W$ have been permuted so that they appear in descending order, and the term corresponding to each row is shown. (Recall the discussion about labelling rows and columns at the end of Section 2.) We also present a few columns of the matrix $H$ for $r=3$.

\begin{table}
\caption{$W$-matrix, $r=2$}
\label{tab:Wr=2}       

\begin{tabular}{lr}
\toprule
function  &  23.0\\
linear  &  11.1 \\
matrix  &  11.0 \\
equation  &  9.9 \\
derivative  &  9.4 \\
calculus  &  7.8 \\
differential  &  6.5 \\
solve  &  6.1\\
problem  &  5.9\\
mathematical  &  5.3 \\
science  &  5.2 \\
compute  &  5.1 \\
variable  &  4.7 \\
applications  &  4.7 \\
integral  &  4.6 \\
\bottomrule
\end{tabular}
\qquad
\begin{tabular}{lr}
\toprule
group  &  17.9\\
theorem  &  15.0 \\
theory  &  6.9 \\
ring  &  6.6 \\
understand  &  4.0\\
structure  &  3.3\\
example  &  3.0 \\
number  &  2.9 \\
isomorphism  &  2.7 \\
concepts  &  2.6 \\
homomorphisms  &  2.5\\
sylow  &  2.3\\
subgroups  &  2.3 \\
quotient  &  2.2 \\
cauchy  &  2.2\\
\bottomrule
\end{tabular}
\end{table}

\begin{table}
\caption{$W$-matrix, $r=3$}
\label{tab:Wr=3}       
\begin{tabular}{lr}
\toprule
group  &  18.1 \\
theorem  &  14.8\\
theory  &  6.8 \\
ring  &  6.7 \\
understand  &  3.9 \\
structure  &  3.4 \\
number  &  3.0 \\
example  &  2.9 \\
isomorphism  &  2.6 \\
homomorphisms  &  2.6\\
concepts  &  2.5\\
sylow  &  2.3 \\
subgroups  &  2.3 \\
applications  &  2.2 \\
quotient  &  2.2 \\
cauchy  &  2.1 \\
time  &  2.1 \\
finite  &  2.1 \\
algebraic  &  2.1 \\
permitting  &  2.0 \\
\bottomrule
\end{tabular}%
\qquad
\begin{tabular}{lr}
\toprule
function  &  27.1 \\
derivative  &  11.1 \\
calculus  &  9.2 \\
equation  &  6.3 \\
differential  &  5.8 \\
integral  &  5.7 \\
problem  &  5.5\\
variable  &  5.3 \\
graph  &  4.8 \\
limit  &  4.6 \\
solve  &  4.5 \\
mathematics  &  4.3 \\
calculate  &  3.9 \\
applications  &  3.8 \\
science  &  3.5 \\
introduction  &  3.5 \\
mathematical  &  3.5 \\
method  &  3.4 \\
polynomial  &  3.4 \\
differentiation  &  3.3\\
\bottomrule
\end{tabular}%
\qquad
\begin{tabular}{lr}
\toprule
matrix  &  19.1 \\
linear  &  14.9 \\
space  &  9.5 \\
vector  &  8.7 \\
algebra  &  7.1 \\
basis  &  6.5 \\
equation  &  6.4 \\
compute  &  5.7\\
system  &  4.9 \\
rank  &  3.5 \\
complex  &  3.5 \\
product  &  3.4 \\
number  &  3.3 \\
mathematical  &  3.3\\
science  &  3.3 \\
solve  &  3.2 \\
dimensional  &  3.0 \\
eigenvalues  &  2.9 \\
set  &  2.8 \\
eigenvectors  &  2.7 \\
\bottomrule
\end{tabular}
\end{table}

\begin{table}
\caption{$W$-matrix, $r=4$}
\label{tab:Wr=4}       
\begin{tabular}[t]{lr}
\toprule
group  &  20.5 \\
theorem  &  11.3 \\
ring  &  6.7 \\
theory  &  4.8 \\
structure  &  3.9 \\
isomorphism  &  3.0 \\
homomorphisms  &  2.9 \\
applications  &  2.7 \\
sylow  &  2.6 \\
subgroups  &  2.6 \\
quotient  &  2.6 \\
algebra  &  2.3 \\
algebraic  &  2.2 \\
time  &  2.1 \\
permitting  &  2.1 \\
finite  &  1.9 \\
lagrange  &  1.8 \\
special  &  1.7 \\
construct  &  1.5 \\
\bottomrule
\end{tabular}%
\qquad
\begin{tabular}[t]{lr}
\toprule
function  &  25.2 \\
derivative  &  11.7 \\
calculus  &  8.4 \\
equation  &  6.7 \\
differential  &  6.3 \\
problem  &  5.8 \\
variable  &  5.5 \\
graph  &  4.9 \\
solve  &  4.8 \\
limit  &  4.3 \\
applications  &  4.2 \\
integral  &  4.1 \\
calculate  &  4.1 \\
polynomial  &  3.7 \\
differentiation  &  3.4 \\
science  &  3.3 \\
mathematics  &  3.2 \\
introduction  &  3.2 \\
inverse  &  3.0 \\
\bottomrule
\end{tabular}%
\qquad
\begin{tabular}[t]{lr}
\toprule
matrix  &  19.3 \\
linear  &  14.8 \\
space  &  9.0 \\
vector  &  8.7 \\
algebra  &  7.1 \\
basis  &  6.5 \\
equation  &  6.4 \\
compute  &  5.7 \\
system  &  4.8 \\
rank  &  3.6 \\
product  &  3.4 \\
solve  &  3.1 \\
science  &  3.1 \\
complex  &  3.0 \\
eigenvalues  &  2.9 \\
dimensional  &  2.9 \\
eigenvectors  &  2.8 \\
number  &  2.7 \\
mathematical  &  2.6 \\
\bottomrule
\end{tabular}%
\qquad
\begin{tabular}[t]{lr}
\toprule
theorem  &  9.0 \\
understand  &  8.7 \\
question  &  6.3 \\
complex  &  6.1 \\
number  &  5.6 \\
example  &  5.6 \\
concepts  &  5.4 \\
mathematical  &  5.2 \\
function  &  5.1 \\
cauchy  &  4.9 \\
theory  &  4.6 \\
integral  &  3.8 \\
demonstrate  &  3.8 \\
correctly  &  3.4 \\
method  &  3.4 \\
series  &  3.3 \\
write  &  3.3 \\
set  &  3.3 \\
sequence  &  3.3 \\
\bottomrule
\end{tabular}
\end{table}

\begin{table}
\caption{$H$-matrix, $r=3$}
\label{tab:Hr=3}    
\begin{tabular}{|@{\hspace{3em}}c@{\hspace{3em}}c@{\hspace{3em}}c@{\hspace{3em}}c|}
\multicolumn{4}{c}{\textbf{Module Names}}\\[1ex]
\multicolumn{1}{c}{\small \shortstack{Numbers and\\ Functions}} &
\small \shortstack{Linear Algebra with \\ Applications to Economics} &
\small \shortstack{Groups, Rings,\\ and Fields} &
\multicolumn{1}{c}{\small \shortstack{Differential Equations\\ via Computer Algebra}}\\ \hline
0.1 & {0.0} & 0.3 & {0.0} \\
0.1 &{0.0} & {0.0} & 0.1 \\
{0.0} & 0.2 & {0.0} & {0.0} \\\hline
\end{tabular}
\end{table}


\clearpage


\section{Discussion}

Following on from the discussion around the small example presented above, the students understand how the frequency of the terms across all documents is computed. This gives an easy and automated way to derive the most and least frequent words, given in Table~\ref{tab:Most Frequent Words}. While in small example shown the most and least frequent words could easily be found by hand, this is impractical when the  matrix becomes large.

\discuss{Consider different columns of the matrix $W$ for $r=2,3$ given in Tables \ref{tab:Wr=2} and \ref{tab:Wr=3}. What do you observe?}

\noindent Already in the case when $r=2$, we can see some regularity in the way the terms are grouped into columns. For example, it makes sense that the terms \emph{function}, \emph{derivative}, \emph{differential}, \emph{integral} appear in the same column. In the second column we see the terms \emph{group}, \emph{ring}, \emph{isomorphism}, \emph{homomorphism}, \emph{sylow}, \emph{subgroups}, \emph{quotient}, \emph{cauchy} appearing together. 

The factorization is perhaps the most informative for the choice $r=3$, so let us take a closer look at this case. The terms in the matrix $W$ are grouped in such a sensible way that we can challenge the students to give them titles. Those students who've read ahead a little may suggest \emph{Abstract Algebra} for the first column, while most should find \emph{Calculus} appropriate for the second, and \emph{Linear Algebra} for the third. Things become a little less clear when we consider the matrix $W$ for $r=4$ (Table~\ref{tab:Wr=4}).

\discuss{What are the advantages and disadvantages of choosing $r$ to be small or big?}

\noindent  While higher values of $r$ will make $\hat{V}$ closer to $V$, they can make it more difficult to interpret the results. An informed choice of $r$, dependent on the needs of the applications, needs to be made.  This problem of ``model complexity'' has been the subject of a great deal of research in Statistics and Machine Learning.

\discuss{In the $r=3$ case, we were able to give titles to columns of the matrix $W$. Those titles could be called ``topics''. The rows of $W$ are indexed by ``words'' and the columns by ``topics''. For the multiplication $WH$ to make sense we need to have the rows of $H$ marked by ``topics''. Let us look at the matrix $H$ given in Table \ref{tab:Hr=3} to see if this makes sense.}

Representative columns given for the matrix $H$ agree with our prediction that the first row corresponds to Abstract Algebra, the second row to Calculus and the third row to Linear Algebra. For example, the course \emph{Linear Algebra with Applications to Economics} has the only nonzero entry in the third row, while the course \emph{Differential Equations via Computer Algebra} has the only nonzero entry in the second row.

\section{Conclusion}

While this is a black box experiment for the students, they are able to appreciate the result and understand the emergence of the topics in an example.
The NMF algorithm yields this topic analysis, helping us appreciate the strengths of the method.
If we want to bring the discussion further, it can be pointed out how this class of algorithm is used to decompose speech and music into phonemes and notes \citep{ASARI-ETAL-2006A, OGRADY-PEARLMUTTER-2008a, Smaragdis-Brown-2003a}, in speech denoising \citep{Wilson-etal-2008a} and recognition \citep{Hurmalainen-2014a}, in chemistry \citep{siy2008matrix} and biomedical sciences \citep{Cancer-NMF, Ortega-Martorell-etal-2012a, Paine-etal-2016a, HIV-NMF}, in the analysis of the cosmic microwave background radiation \citep{CARDOSO03}, etc.

The example presented above can be adapted for classroom needs in various ways. 
An aspect not discussed here is the potential to turn some of the above ideas into student projects.
We are aware that the computational aspects of this may be a big stumbling block, so we are developing a web-based tool to make it easy for students to analyse a set of document in this way. We see a potential for interdisciplinary projects, where students are charged with the task of analyzing a large body of documents on a particular subject, and use linear algebra to reach some conclusions. 

In collaboration with Miao Wei,\footnote{Dept of Computer Science, Maynooth University, Ireland, \email{davidweimiao@gmail.com}.
} we have created an end-to-end interactive browser-based implementation of the processing pipeline discussed above (taking documents as input and processing them through stemming, the construction of a term-by-document matrix, NMF, and visualization of the resulting factor matrices), which is being made available online.\footnote{\url{http://barak.pearlmutter.net/demo/NMF/}}


\bibliographystyle{plainnat}    
\bibliography{References}       
\end{document}